\documentclass[reprint,aps]{revtex4-2}
\usepackage[english]{babel}
\usepackage{amsmath,amssymb,enumerate,mathrsfs,accents,xcolor,tabularx}

\newcommand{\assign}{:=}
\newcommand{\backassign}{=:}
\newcommand{\cdummy}{\cdot}
\newcommand{\mathD}{\mathrm{D}}
\newcommand{\mathd}{\mathrm{d}}
\newcommand{\tmcolor}[2]{{\color{#1}{#2}}}
\newcommand{\tmem}[1]{{\em #1\/}}
\newcommand{\tmop}[1]{\ensuremath{\operatorname{#1}}}
\newcommand{\tmtextbf}[1]{\textbf{{{#1}}}}
\newcommand{\tmtextit}[1]{\textit{{{#1}}}}
\newenvironment{enumeratealpha}{\begin{enumerate}[a{\textup{)}}] }{\end{enumerate}}
\newtheorem{definition}{Definition}
\newtheorem{proposition}{Proposition}

\newcommand{\bbE}{\mathbb{E}}

\begin{document}

\title{Wilson--It{\^o} diffusions}

\author{Ismael Bailleul}
\email[Email: ]{ismael.bailleul@univ-brest.fr}
\affiliation{Univ Brest, CNRS UMR 6205,\\
Laboratoire de Math{\'e}matiques de Bretagne Atlantique, France}

\author{Ilya Chevyrev}
\email[Email: ]{ichevyrev@gmail.com}
\affiliation{School of Mathematics, University of Edinburgh, United Kingdom}

\author{Massimiliano Gubinelli}
\email[Email: ]{gubinelli@maths.ox.ac.uk}
\affiliation{Mathematical Institute, University of Oxford, United Kingdom}


\begin{abstract}
  We introduce Wilson--It{\^o} diffusions, a class of random fields on
  $\mathbb{R}^d$ that change continuously along a scale parameter via a
  Markovian dynamics with \tmtextit{local} coefficients. Described via
  forward-backward stochastic differential equations, their observables
  naturally form a pre-factorization algebra \tmtextit{{\`a} la}
  Costello--Gwilliam. We argue that this is a new non-perturbative
  quantization method applicable also to gauge theories and independent of a
  path-integral formulation. Whenever a path-integral is available, this
  approach reproduces the setting of Wilson--Polchinski flow equations.
\end{abstract}

{\maketitle}

\section{Introduction}

The Kadanoff--Wilson point of view on quantum field theory is a central idea
in modern physics. It dictates that one should develop the fluctuation of
quantum or statistical fields along a scale decomposition associated to the
tunable precision of the observation device. In this framework a random
Euclidean field is described by a stochastic process $(\varphi_a)_{a \geqslant
0}$ where $a$ parametrizes a scale of observation with characteristic length
$1 / a$ and $\varphi_a$ is the corresponding observation of the field, which
must be thought of as containing fluctuations with spatial scales $\gtrsim 1 /
a$. It is natural to assume that, as we gradually increase the resolution of
our measuring devices, the resulting measurements vary continuously. Therefore
we will postulate that the stochastic process $\varphi_a$ is a
{\tmem{pathwise}} continuous function of $a$. It is also reasonable to assume
that $\varphi_a$ contains all the information gathered by observations with
lesser precisions, which implies that it is a Markov process along the scale
parameter, and we assume further that $\varphi_a \rightarrow \varphi_{\infty}$
as $a \rightarrow \infty$. We denote by $\mathscr{F} = (\mathscr{F}_a)_{a
\geqslant 0}$ the filtration generated by $(\varphi_a)_{a \geqslant 0}$ and by
$\bbE_b [\hspace{0.17em} \cdummy \hspace{0.17em}]$ the operator of conditional
expectation given $\mathscr{F}_b$, for any $b \geqslant 0$, with
$\mathbb{E}=\mathbb{E}_0$ the full expectation. We assume that $\varphi_0 = 0$
whenever the fields take value in a vector space; otherwise we take
$\varphi_0$ to be a fixed default (classical) configuration (e.g. a background
field). We say that a functional $\varphi \mapsto F (\varphi)$ of the field
$\varphi$ is supported in $U \subseteq \mathbb{R}^d$ if $F (\varphi) = F
(\psi)$ for any field $\psi$ such that $\varphi = \psi$ in $U$. For
$\varepsilon > 0$ let $U_{\varepsilon} \assign \{ x \in \mathbb{R}^d : d (x,
U) < \varepsilon \}$ be the $\varepsilon$-enlargement of $U$.

\begin{definition}
  \label{def:obs}An \tmtextbf{observable} $O$ is a stochastic process
  $(O_a)_{a \geqslant 0}$ which is an $\mathscr{F}$-martingale, i.e.
  \[ O_b = \bbE_b [O_a], \qquad (0 \leqslant b \leqslant a < \infty) . \]
  An \tmtextbf{observable} $O$ is said to be \tmtextbf{supported on} $U$ if
  there exists a family of functionals $(\ring{O}_a (\cdot))_{a \geqslant 0}$
  each supported on the $a^{- 1}$-enlargement $U_{a^{- 1}}$ and such that
  $R_a^O \assign O_a - \ring{O}_a (\varphi_a) \rightarrow 0$ as $a \rightarrow
  \infty$. We call such $\ring{O}$ a \tmtextbf{germ for $O$}. A
  \tmtextbf{local observable (field)} is field of observables $x \mapsto O
  (x)$ such that the observable $O (x)$ is supported on the singleton $\{x\}$
  for all $x \in \mathbb{R}^d$.
\end{definition}

The definition of the $a^{- 1}$-enlargement is a matter of convention. An
observable $O$ that has a germ
\[ \ring{O}_a (\varphi_a) (x) =\mathcal{O}_a (\varphi_a (x), \ldots, \nabla^k
   \varphi_a (x)), \]
for some finite $k$ and some function $\mathcal{O}_a$, is local. Germs serve
to parametrize the space of local observables via concrete functionals of the
fields.

One owes to \tmtextbf{\tmtextit{It{\^o}}} the fundamental insight that
continuous diffusions can be constructed via stochastic differential equations
{\cite{stroockMarkovProcessesIto2003}}. As such, since we postulated that the
scale description $(\varphi_a)_a$ of a random field is a \ continuous Markov
process, it is completely described via its infinitesimal rate of change
\[ \mathd \varphi_a = \varphi_{a + \mathd a} - \varphi_a = B^{\varphi}_a
   \mathd a + \mathd M^{\varphi}_a \]
which comes in two parts: a \tmtextit{drift} $B^{\varphi}_a \mathd a$ and a
\tmtextit{martingale} part $\mathd M^{\varphi}_a$. The drift component
$B^{\varphi}_a$ is the part of the rate of change which can be predicted from
the large scale observations, while the ``innovation'' process $M$ is a
martingale w.r.t. this filtration, it models the additional information gained
by augmenting the resolution. Recall a continuous martingale is described via
its quadratic variation process $\langle M
\rangle_a$~{\cite{Rogers_Williams_2000}}.

Enters \tmtextbf{\tmtextit{Wilson}}. We inject in this standard framework the
spatial random structure of Euclidean quantum fields by imposing that the
dynamics is completely described by `\tmtextit{local}' coefficients. We do so
by postulating that there exists a local observable field $(f_a)_{a \geqslant
0}$ which models the \tmtextbf{\tmtextit{microscopic force}} and that,
inspired by the Kadanoff--Wilson block averaging procedure, the drift
$B^{\varphi}_a$ at scale $a$ is determined by an averaging of the force field
$f_a$ over a region of size $1 / a$. Therefore we introduce a ``block
averaging'' operator $C_a = C_a (\varphi_a)$, which could be random and depend
on $\varphi_a$. We require $C_a$ to have support on a ball of radius $1 / a$,
that $a \mapsto C_a (\psi)$ varies smoothly for all $\psi$, that $C_a$ be
symmetric and positive definite and such that $C_{\infty} = 1$. For example,
one can take
\begin{equation}
  (C_a h) (x) = \int a^d \chi ((x - y) a) h (y) \mathd y,
  \label{eq:main-averaging}
\end{equation}
where $x \in \mathbb{R}^d$ and $\chi$ is a smooth, radially symmetric,
positive definite function of unit integral. In particular note that $C_0 = 0$
on all sufficiently nice functions $h$. Since the drift will be integrated
along the scales, the local averaging has to be done in such a way as to not
over-count the contributions of the microscopic force. We denote by $\dot{C}_a
\assign \partial_a C_a$ the scale-derivative of this averaging and let
$B^{\varphi}_a \assign \dot{C}_a f_a$. This fixes the previsible part of the
stochastic dynamics as a function of the microscopic force. To complete our
description we need to specify also the quadratic variation of the martingale
part. Using the same principles we can assume
\[ \mathd \langle M \rangle_a \assign \dot{C}^{1 / 2}_a \sigma^2_a 
   \dot{C}^{1 / 2}_a \mathd a \]
for a microscopic positive ``diffusivity'' $(\sigma^2_a)_a$ which is a local
observable field of positive scalars. In this way the local diffusion on
scales is determined by the datum of two local observable fields $(f_a,
\sigma^2_a)_a$ and a family of averaging operators $(C_a)_a$. In particular,
by standard results there exists, possibly on an extended probability space, a
cylindrical $\mathscr{F}$-Brownian motion
\[ \mathbb{E} [W_b (x) W_a (y)] = (b \wedge a) \delta (x - y), \]
such that
\[ M_a = \int_0^a \dot{C}^{1 / 2}_b \sigma_b \mathd W_b \]
where $\sigma_b \assign (\sigma^2_b)^{1 / 2}_b$.

\begin{definition}
  \label{DefnLocalEQFT}A \tmtextbf{Wilson--It{\^o} diffusion} is a continuous
  stochastic process $(\varphi_a)_{a \geqslant 0}$ taking values in the set of
  smooth functions on $\mathbb{R}^d$ with the following properties.
  \begin{enumeratealpha}
    \item \tmtextbf{Dynamics}. There is an \tmtextit{effective force}
    $(f_a)_{a \geqslant 0}$ and \tmtextit{an effective diffusivity}
    $(\sigma_a^2)_{a \geqslant 0}$ such that $(\varphi_a)_{a \geqslant 0}$ is
    a Markovian It{\^o} diffusion
    \begin{equation}
      \label{EqWIDE} \mathd \varphi_a = \dot{C}_a f_a \mathd a + \dot{C}_a^{1
      / 2} \sigma_a \mathd W_a .
    \end{equation}
    \item \tmtextbf{Locality}. The effective force $f$ and the effective
    diffusivity $\sigma^2$ are local observable fields.
  \end{enumeratealpha}
  We call equation~\eqref{EqWIDE} a \tmtextbf{Wilson--It{\^o} differential
  equation} (WIDE). A \tmtextbf{Wilson--It{\^o} field} is the random field
  $\varphi_{\infty}$ obtained as the terminal value of a Wilson--It{\^o}
  diffusion $(\varphi_a)_{a \geqslant 0}$.
\end{definition}

The main goal of this paper is to propose the hypothesis that Euclidean
quantum field theories can be identified with Wilson--It{\^o} fields. This
provides a new framework, independent of the path-integral formalism, to study
Euclidean quantum fields.
\begin{enumeratealpha}
  \item This description emerges from simple and natural assumptions and
  covers in principle much more than those theories that can be reached
  perturbatively from a Gaussian functional integral.
  
  \item The continuity of the process $(\varphi_a)_a$ with respect to $a$
  gives it the structure of an It{\^o} diffusion. When $\sigma_a$ and $C_a$
  are deterministic (hence $\sigma_a$ is constant by the martingale property
  of $\sigma_a^2$), and without any assumption of a perturbative regime, the
  Wilson--It{\^o} random field $\varphi_{\infty}$ comes with an associated
  Gaussian field
  \[ X^C_a \assign \int_0^a \dot{C}^{1 / 2}_b \sigma_b \mathd W_b, \]
  and a coupling to it. In that case, $X^C_{\infty}$ is a white noise
  \[ \mathbb{E} [X^C_{\infty} (x) X^C_{\infty} (y)] = \sigma^2 C_{\infty} (x -
     y) = \sigma^2 \delta (x - y) . \]
  Therefore one expects that
  \[ \varphi_{\infty} = X_{\infty}^C + \int_0^{\infty} \dot{C}_a f_a \mathd a
  \]
  is, in general, only a distribution of very low regularity.
  
  \item Eq.~\eqref{EqWIDE} makes sense for fields defined on and/or taking
  values in manifolds or vector bundles and for which the path-integral
  formalism is less clear to apply. Our formalism is non-perturbative and
  trades the use of functional integrals against It{\^o} calculus.
  
  \item Let $A = A (a)$ be a possibly random, adapted, increasing change of
  scale such that $A (0) = 0$ and $A (\infty) = \infty$, and let $\mathd
  \tilde{W}_a \assign A' (a)^{- 1 / 2} \mathd W_{A (a)}$, $\tilde{C}_a = C_{A
  (a)}$, $\tilde{f}_a \assign f_{A (a)}$ and $\tilde{\sigma}_a \assign
  \sigma_{A (a)}$. Then we have
  \[ \mathd (\varphi_{A (a)}) = \partial_a \tilde{C}_a  \tilde{f}_a \mathd a +
     (\partial_a \tilde{C}_a)^{1 / 2}  \tilde{\sigma}_a \mathd \tilde{W}_a, \]
  which shows that Wilson--It{\^o} diffusions are covariant wrt. random
  changes of spatial scales. In particular, this justifies that the diffusion
  has to be averaged with $\dot{C}_b^{1 / 2}$. Note that observables are also
  covariant: If $(O_a)_a$ is an observable for $(\varphi_a)_a$ then $(O_{A
  (a)})_a$ is an observable for the diffusion $(\varphi_{A (a)})_a$.
  
  \item Only the law of the terminal value $\varphi_{\infty}$ and the averages
  of observables are the physical content of a Wilson--It{\^o} diffusion.
\end{enumeratealpha}
\tmtextbf{Dyson--Schwinger equations and martingale problems} -- The law of
the process $(\varphi_a)_{a \geqslant 0}$ is determined by a
\tmtextbf{\tmtextit{martingale problem}}: for all sufficiently nice
scale-dependent test functions $F (a, \psi)$, the process
\[ M^F_a \assign F (a, \varphi_a) - \int_0^a (\mathscr{L}_b F)  (b, \varphi_b)
   \mathd b \]
is an $\mathscr{F}$-martingale where
\begin{equation}
  \begin{array}{lll}
    \mathscr{L}_b F (b, \varphi) & \assign & \partial_b F (b, \varphi) +
    \mathD F (b, \varphi)  \dot{C}_b f_b \\
    &  & + \frac{1}{2} \tmop{Tr} [\dot{C}_b^{1 / 2} \sigma_b^2 \dot{C}_b^{1 /
    2} \mathD^2 F (b, \varphi)]
  \end{array} \label{eq:generator}
\end{equation}
is the \tmtextbf{generator} of the Wilson--It{\^o} diffusion. Here we denote
by $\mathD \assign \delta / \delta \varphi$ the gradient of a functional
$\varphi \mapsto F (\varphi)$ with respect to the field $\varphi$ taken wrt.
the $L^2$ norm and by $\tmop{Tr}$ the trace operator.

\section{Forward-backward stochastic differential equations}

Some observables can be obtained via the conditional expectation of a function
of $\varphi_{\infty}$, i.e. let $O_a =\mathbb{E}_a [F (\varphi_{\infty})]$.
This kind of martingale is called \tmtextit{closed}, i.e. it admits a
``terminal value'' $O_{\infty} = F (\varphi_{\infty})$ from which it can be
reconstructed. In the case when $F$ is a linear functional we have
\[ O_a =\mathbb{E}_a [F (\varphi_{\infty})] = F (\mathbb{E}_a
   [\varphi_{\infty}]), \]
assuming suitable integrability conditions here and below. Using the
Wilson--It{\^o} dynamics, and since $(f_a)_a$ is a martingale, we have
\begin{equation}
  \begin{array}{lll}
    \mathbb{E}_a [\varphi_{\infty}] & = & \varphi_a +\mathbb{E}_a
    \int_a^{\infty} \dot{C}_b f_b \mathd b = \varphi_a + \int_a^{\infty}
    \dot{C}_b f_a \mathd b\\
    & = & \varphi_a + C_{\infty, a} f_a
  \end{array} \label{eq:lin-comp}
\end{equation}
where \ $C_{\infty, a} \assign C_{\infty} - C_a$. In general we are not
allowed to form non-linear \tmtextit{local} functions $F (\varphi_{\infty})$
of the distribution $\varphi_{\infty}$: the locality condition for an
observable is non-trivial and usually requires renormalization. As a
consequence, we do not expect non-linear local observables to be closed
martingales. Ignoring for the moment this difficulty, we note that when $F
(\varphi_{\infty})$ is well-defined the observables $(O_a)_a$ are closed
martingales and they satisfy a \tmtextit{backward} \tmtextit{stochastic
differential equations}
(BSDEs)~{\cite{maForwardBackwardStochasticDifferential2007}}
\[ \mathd O_a = Z^O_a \mathd W_a, \]
for a pair of adapted processes $(O_a, Z^O_a)$ with \tmtextit{terminal}
condition $O_{\infty} = F (\varphi_{\infty})$.

A procedure to construct local observables starts with some approximate local
observable given by a function $\ring{O}^{a_0} (\varphi_{\infty})$ localized
at scale $a_0^{- 1}$ and setting
\[ O^{a_0}_a \assign \bbE_a [\ring{O}^{a_0} (\varphi_{\infty})] . \]
One can then study the convergence of the family of the \tmtextit{non-local}
observables $(O^{a_0}_a)_{0 \leqslant a \leqslant a_0}$ as $a_0 \rightarrow
\infty$. Provided the functions $\ring{O}^{a_0} (\varphi_{\infty})$ contains
appropriate (diverging) renormalizations one is able to show that the
observables $O^{a_0}$ converge to a local observable as $a_0 \rightarrow
\infty$.

This approach leads naturally to the analysis of a general class of
\tmtextit{forward-backward stochastic differential equation} (FBSDEs) of the
form
\begin{equation}
  \mathd \phi_a = \dot{C}_a \tilde{f}_a \mathd a + \dot{C}_a^{1 / 2}
  \tilde{\sigma}_a \mathd W_a \label{eq:approx-sde}
\end{equation}
for some approximately local functionals
\[ \tilde{f}_a = \bbE_a [F (\phi_{\infty})], \quad \tilde{\sigma}_a = \bbE_a
   [\Sigma^2 (\phi_{\infty})]^{1 / 2} . \]
Recall $\dot{C}_a$ may depend on $\phi_a$ in a general setting. Note also that
we use the letter $\phi$ for these ``approximate'' dynamics while we use
$\varphi$ for a local dynamics. These FBSDE are not proper Wilson--It{\^o}
diffusions, since their coefficients are non-local. BSDEs and FBSDEs are well
studied in the mathematical literature (see
e.g.~{\cite{maForwardBackwardStochasticDifferential2007}}) and we dispose also
of numerical methods to approximate their solutions. In relation to the
numerical aspects, note that the formalism allows us to replace $\mathbb{R}^d$
by a finite discrete lattice $(\varepsilon \mathbb{Z} \cap [- L, L])^d$ with
mesh $\varepsilon$.

\subsection{Linear-like force}

To give a first example consider the case of an approximate force $F$ with a
linear component and a constant diffusivity
\[ F (\phi_{\infty}) = \tmcolor{red}{\tmcolor{black}{\alpha}} (- A
   \phi_{\infty} + h (\phi_{\infty})), \qquad \Sigma (\phi_{\infty}) =
   \tmcolor{red}{\tmcolor{black}{\alpha^{1 / 2}}}, \]
for some positive constant $\alpha$, some positive linear operator $A$ and
some additional force component $h (\phi_{\infty})$. For a local operator $A$
the linear functional $A \phi_{\infty}$ is always well-defined in the space of
distributions, so it defines a local force field. We assume here that the
averaging operators $C_a$ are deterministic, field-independent and commute
with $A$. Think of the case where $A = m^2 - \Delta$, $h = 0$ and $C_a$ is
given by Eq.~\eqref{eq:main-averaging}. Using~\eqref{eq:lin-comp} we have
\[ \mathbb{E}_a [\phi_{\infty}] = \tmcolor{black}{\phi_a - \alpha C_{\infty,
   a} A\mathbb{E}_a [\phi_{\infty}] + \alpha C_{\infty, a} \mathbb{E}_a [h
   (\phi_{\infty})] } \]
Solving for $\mathbb{E}_a [\phi_{\infty}]$ and letting
\[ \psi_a \assign (1 + \alpha C_{\infty, a} A)^{- 1} \phi_a, \]
we have $\psi_{\infty} = \phi_{\infty}$ and
\begin{equation}
  \mathd \psi_a = \dot{Q}_a \bbE_a [h (\psi_{\infty})] \tmcolor{black}{\mathd
  a + \dot{Q}_a^{1 / 2}} \mathd W_a, \label{eq:Q-eq}
\end{equation}
where $\dot{Q}_a \assign \tmcolor{black}{\partial_a}  (A^{- 1} (1 + \alpha
C_{\infty, a} A)^{- 1})$. This computation shows that the linear component in
the force can always be integrated and gives rise to a modified FBSDE where
the local averaging operator $\dot{C}_a$ has been replaced by the operator
$\dot{Q}_a$. Note in particular that the Gaussian field
\begin{equation}
  X_a^Q \assign \int_0^a \dot{Q}_c^{1 / 2} \mathd W_c, \label{eq:GFF}
\end{equation}
has covariance $Q_a \tmcolor{black}{- Q_0}$ and that $X_{\infty}^Q$ has
covariance $\alpha (1 + \alpha A)^{- 1}$. One gets back the operator $A^{- 1}$
in the large $\alpha$ limit, in which case $X_{\infty}^Q$ is a massive GFF in
the model situation. We can invert the transformation and go from a FBSDE of
the form~\eqref{eq:Q-eq} back to the FBSDE~\eqref{eq:approx-sde} noting that
\[ Q_{\infty, a}^{- 1} - \alpha^{- 1} C_{\infty, a}^{- 1} = A \]
for all $a \geqslant 0$.

\subsection{Gradient diffusions}

Assume further that $A$ is a symmetric operator and the additional force
component $h$ is given by an effective UV-regularized potential $V_{\infty}$
as
\[ h (\psi_{\infty}) = - \mathD V_{\infty} (\psi_{\infty}), \]
and let
\[ V_a (\varphi) \assign - \log \mathbb{E} [e^{- V_{\infty} (\varphi +
   X_{\infty}^Q - X_a^Q)}] . \]
Then $(V_a)_a$ is the solution of the Polchinski flow
equation~{\cite{polchinskiRenormalizationEffectiveLagrangians1984}}
\begin{equation}
  \partial_a V_a - \frac{1}{2} \mathD V_a \dot{Q}_a \mathD V_a + \frac{1}{2}
  \dot{Q}_a \mathD^2 V_a = 0 \label{eq:polchinski}
\end{equation}
with terminal condition $V_{\infty}$ for $a \rightarrow \infty$. Letting $Z_a
\assign \exp (V_a (\psi_a) - V_0 (0))$ and using It{\^o} formula, we have that
$Z_a$ is a positive martingale with $Z_0 = 1$ and moreover that, under the
probability measure $\mathbb{Q}$ defined on $\mathscr{F} _a$ by $\mathd
\mathbb{Q} \assign Z_a \mathd \mathbb{P}$, the process $\psi_a$ is a
martingale with deterministic quadratic variation $\mathd \langle \psi
\rangle_a = \dot{Q}_a \mathd a$, so $\psi_{\infty}$ is under $\mathbb{Q}$ a
GFF. Note also that it follows from~\eqref{eq:polchinski} that $F_a (\psi_a)
\assign - \mathD V_a (\psi_a)$ is a martingale under $\mathbb{Q}$. Recalling
Eq.~\eqref{eq:GFF}, we conclude that
\begin{equation}
  \begin{array}{lll}
    \mathbb{E}_{\mathbb{P}} [G (\psi_a)] & = & \mathbb{E}_{\mathbb{Q}} [G
    (\psi_a) e^{V_0 (0) - V_a (\psi_a)}]\\
    & = & \mathbb{E}_{\mathbb{P}} [G (X_a^Q) e^{V_0 (0) - V_a (X_a^Q)}]\\
    & = & \frac{\mathbb{E}_{\mathbb{P}} [G (X_a^Q) e^{- V_{\infty}
    (X_{\infty}^Q)}]}{\mathbb{E}_{\mathbb{P}} [e^{- V_{\infty}
    (X_{\infty}^Q)}]}
  \end{array} \label{eq:gibbs}
\end{equation}
for any function $G$ and any $a \geqslant 0$. Here we used that $e^{- V_a
(X^Q_a)}$ is a martingale. This implies that the law $\nu_{\infty}$ of the
random field $\psi_{\infty} = \phi_{\infty}$ is given by
\begin{equation}
  \nu_{\infty} (\mathd \psi) = \frac{e^{- V_{\infty} (\psi)} \mu^{Q_{\infty}}
  (\mathd \psi)}{\int e^{- V_{\infty} (\psi)} \mu^{Q_{\infty}} (\mathd \psi)}
  \label{eq:gibbs-two}
\end{equation}
where $\mu^{Q_{\infty}}$ denotes the Gaussian of covariance $Q_{\infty, 0}$.

\

This shows that the class of Wilson--It{\^o} fields comprise as a particular
case the Euclidean quantum fields~\eqref{eq:gibbs-two} constructed as
perturbations of a Gaussian field. They are obtained by solving FBSDEs of the
form
\begin{equation}
  \mathd \psi_a = - \dot{Q}_a \bbE_a [\mathD V_{\infty} (\psi_{\infty})]
  \mathd a + \dot{Q}_a^{1 / 2} \mathd W_a . \label{eq:fbsde}
\end{equation}
which we call \tmtextbf{\tmtextit{Polchinski FBSDEs}}, since they describe the
Polchinski
semigroup~{\cite{bauerschmidtLogSobolevInequality2021,bauerschmidtStochasticDynamicsPolchinski2023}}.
Even in this potential framework where there is a formal link with the
well-known Polchinski flow equation, the FBSDE~\eqref{eq:fbsde} does not
require the a priori knowledge of the solution of Eq.~\eqref{eq:polchinski}.
Indeed we only used it to derive the formal connection while in
Eq.~\eqref{eq:fbsde} only the boundary condition $\mathD V_{\infty}$ is
needed.\quad In general this quantity needs an ultraviolet regularization that
has to be tuned in order for the FBSDE~\eqref{eq:fbsde} to reach a
well-defined local limit as the regularization is removed.

\subsection{One-dimensional gradient diffusions}

Take $d = 1$ and let
\[ F (\varphi_{\infty}) (x) = - (- \Delta) \varphi_{\infty} (x) - \rho (x) v'
   (\varphi_{\infty} (x)) \]
$x \in \mathbb{R}$, where $v'$ is the derivative of a function $v : \mathbb{R}
\rightarrow \mathbb{R}$ smooth and bounded from below and $\rho : \mathbb{R}
\rightarrow \mathbb{R}_+$ is a compactly supported function of the space
variable (an IR cutoff). The formulation~\eqref{eq:Q-eq} with $A = - \Delta$
reads
\begin{equation}
  \mathd \psi_a = - \dot{Q}_a \rho \mathbb{E}_a v' (\psi_{\infty}) \mathd a +
  \dot{Q}_a^{1 / 2} \mathd W_a . \label{eq:mild-1d}
\end{equation}
From this equation we see that $\psi_{\infty}$ is actually comparable in
regularity to the $d = 1$ massive GFF $X_{\infty}$ with covariance $(1 + A)^{-
1}$, i.e.
\[ \mathbb{E} [X_{\infty} (x) X_{\infty} (y)] = \frac{1}{2} e^{- | x - y |},
   \qquad x, y \in \mathbb{R}, \]
which is a H{\"o}lder continuous function for which the mild
formulation~\eqref{eq:mild-1d} makes sense. The residual force is given by the
well-defined potential $V_{\infty} (\varphi) = \int_{\mathbb{R}} \rho (x) v
(\varphi (x)) \mathd x$. By~\eqref{eq:gibbs} we have
\begin{equation}
  \mathbb{E}_{\mathbb{P}} [G (\psi_{\infty})] = \frac{\mathbb{E}_{\mathbb{P}}
  [G (X_{\infty}) e^{- V_{\infty} (X_{\infty})}]}{\mathbb{E}_{\mathbb{P}}
  [e^{- V_{\infty} (X_{\infty})}]} \label{eq:one-dim-path-integral}
\end{equation}
and this shows that $\psi_{\infty}$ is a (time-inhomogeneous) Markovian
diffusion in the space parameter since $X_{\infty}$ can be described also as
the Ornstein--Uhlenbeck process, solution of the stochastic differential
equation
\[ \mathd X_{\infty} (x) = - X_{\infty} (x) \mathd x + \mathd B_x, \]
where $(B_x)_{x \in \mathbb{R}}$ is a two-sided Brownian motion. Note that
this is a Markov process in the space variable. Moreover if we take $\rho
\rightarrow 1$, under suitable assumptions, the dynamics will converge to a
time-homogeneous space-Markovian diffusion described by the WIDE (recall that
$\varphi_{\infty} = \psi_{\infty}$)
\[ \mathd \varphi_a = \dot{C}_a \mathbb{E}_a [- (- \Delta) \varphi_{\infty} -
   v' (\varphi_{\infty})] \mathd a + \dot{C}_a^{1 / 2} \mathd W_a, \]
along the scale variable $a \geqslant 0$, and, at the same time, it is also
described by the stochastic differential equation
\[ \mathd \varphi_{\infty} (x) = \partial_{\phi} \log \Psi (\varphi_{\infty}
   (x)) \mathd x + \mathd B_x, \]
along the space variable $x \in \mathbb{R}$, where $\Psi : \mathbb{R}
\rightarrow \mathbb{R}_+$ is the ground state of the Hamiltonian with
potential $\phi^2 / 2 + v (\phi)$ (see
e.g.~{\cite{lorincziFeynmanKactypeTheoremsGibbs2011}}).

In absence of the IR cutoff $\rho$, the random field $\psi_{\infty}$ is not
absolutely continuous wrt. $X_{\infty}$ and therefore the path-integral
formulation~\eqref{eq:one-dim-path-integral} loses its meaning while the WIDE
formulation remains valid.

This example shows that a local gradient Wilson--It{\^o} dynamics gives rise
to a Markov process in the space variable. It is natural to conjecture that
this is a general feature of (a wide class of) higher-dimensional
Wilson--It{\^o} diffusions.

\subsection{Variational formulation}

Eq.~\eqref{eq:fbsde} can be interpreted as the Euler--Lagrange equation for a
stochastic control problem. To derive this problem we test the process
\[ u_a \assign - \dot{Q}_a^{1 / 2} \bbE_a [\mathD V_{\infty} (\psi_{\infty})]
\]
with an adapted test field $(v_a)_a$ and integrate both in scale and in the
probability space to get the equation
\begin{equation}
  \mathbb{E} \left[ \int_0^{\infty} \langle v_a, u_a \rangle \mathd a +
  \left\langle \int_0^{\infty} \dot{Q}_a^{1 / 2} v_a \mathd a, \mathD
  V_{\infty} (\psi_{\infty}) \right\rangle \right] = 0 \label{eq:weak}
\end{equation}
Eq.~\eqref{eq:weak} is the first order condition for the problem of minimizing
the functional
\begin{equation}
  \Psi (u) \assign \mathbb{E} \left[ V_{\infty} (\psi_{\infty}^u) +
  \frac{1}{2} \int_0^{\infty} \langle u_a, u_a \rangle \mathd a \right]
  \label{eq:var}
\end{equation}
over all adapted controls $(u_a)_{a \geqslant 0}$, where
\[ \psi_a^u \assign \int_0^a \dot{Q}_b^{1 / 2} u_b \mathd b + \int_0^a
   \dot{Q}_b^{1 / 2} \mathd W_b, \]
is the controlled process.

\subsection{Rigorous results}

Variants of Eq.~\eqref{eq:fbsde}, Eq.~\eqref{eq:weak} or of the variational
problem in Eq.~\eqref{eq:var} have been used to construct several Euclidean
quantum fields including the $\Phi^4_2$ and $\Phi^4_3$
models~{\cite{barashkov2021variational,barashkov2020variational,barashkov2021phi43}},
the H{\o}egh-Krohn model~{\cite{barashkov2021variational}}, the Sine--Gordon
model~{\cite{barashkovStochasticControlApproach2022}} and certain subcritical
Euclidean fermionic field theories~{\cite{de_vecchi22}}. This shows that our
approach is intrinsically non-perturbative and can be made rigorous. We invite
the reader to compare this situation with the non-trivial mathematical
difficulties of the path-integral formalism without cutoffs. As an example,
take the $\Phi^4_3$ Euclidean quantum field on a torus, constructed
in~{\cite{barashkov2020variational}} via a slightly different version of the
variational formulation~\eqref{eq:var}. It is known, and proven
in~{\cite{barashkov2021phi43}}, that the $\Phi^4_3$ measure is not absolutely
continuous with respect to the Gaussian free field, so there cannot be a
rigorous path-integral for it. Similarly, in~{\cite{barashkov2021variational}}
it is shown that some variants of Eq.~\eqref{eq:weak} provide effective tools
to study the infinite volume limit of the $\Phi^4_2$ and of the $\exp
(\Phi)_2$ Euclidean fields.

\section{Properties of Wilson--It{\^o} diffusions}\label{SectionConstruction}

\subsection{Coherent germs}\label{SectionCoherence}

Let $\ring{O}$ be a germ for an observable $O$. For $0 \leqslant b \leqslant
a$, one has by It{\^o}'s formula
\[ \begin{array}{lll}
     O_b & = & \ring{O}_b (\varphi_b) + R^O_b = \bbE_b  [\ring{O}_a
     (\varphi_a) + R^O_a]\\
     & = & \ring{O}_b (\varphi_b) + \bbE_b \left[ \int_b^a \mathscr{L}_c
     \ring{O}_c (\varphi_c) \mathd c \right] + \bbE_b [R^O_a]
   \end{array} \]
where $\mathscr{L}_c$ is the generator~\eqref{eq:generator} of the
Wilson--It{\^o} diffusion. The assumption that $R^O$ goes to $0$ in a strong
enough sense gives
\begin{equation}
  R^O_b = \bbE_b \left[ \int_b^{\infty} \hspace{-0.1cm} \mathscr{L}_c
  \ring{O}_c (\varphi_c) \mathd c \right] . \label{eq:remO}
\end{equation}
Therefore $R^O$, and hence $O$ itself, is completely determined by the germ
$\ring{O}$ provided the integral
\[ \int_b^{\infty} \left\|| \bbE_b [\mathscr{L}_c \ring{O}_c (\varphi_c)]
   \right\|| \mathd c \]
converges absolutely. A germ which has this property is called a
\tmtextit{\tmtextbf{coherent germ}}; its associated observable is determined
from it. Note that Eq.~\eqref{eq:remO} for the remainder is equivalent to the
BSDE
\begin{equation}
  \mathd R^O_a = - \hspace{-0.1cm} \mathscr{L}_a \ring{O}_a (\varphi_a) \mathd
  a - Z^O_a \mathd W_a, \quad R^O_{\infty} = 0 \label{eq:bsde-O}
\end{equation}
for the pair of adapted processes $(R^O_a, Z^O_a)_a$.

\

The effective force itself is an observable. For simplicity we assume that
the diffusivity is taken constant $\sigma_a^2 = 1$, similar considerations
otherwise apply to it. Assume that the force has a germ $\ring{f}_a
(\varphi_a)$ and a remainder $R^f$ which then, due to~\eqref{eq:remO},
satisfies
\begin{equation}
  \begin{array}{lll}
    R^f_b & = & \int_b^{\infty} \bbE_b  [R^f_c  \dot{C}_c \mathD \ring{f}_c
    (\varphi_c)] \mathd c\\
    &  & + \int_b^{\infty} \bbE_b [\ring{\mathscr{L}}_c  \ring{f}_c
    (\varphi_c)] \mathd c
  \end{array} \label{EqRemainderForceField}
\end{equation}
where we introduced the operator
\[ \ring{\mathscr{L}_c} \assign \partial_c + \ring{f}_c  \dot{C}_c \mathD +
   \frac{1}{2} \tmop{Tr} \dot{C}_c \mathD^2  . \]
Similarly to the general case of Eq.~\eqref{eq:bsde-O},
Eq.~\eqref{EqRemainderForceField} gives rise to a BSDEs for the pair $(R^f_a,
Z^f_a)_a$ which reads
\begin{equation}
  \begin{array}{lll}
    \mathd R^f_a & = & - \ring{\mathscr{L}}_a \ring{f}_a (\varphi_a) \mathd a
    - R^f_a  \dot{C}_a \mathD \ring{f}_a (\varphi_a) \mathd a\\
    &  & - Z^f_a \mathd W_a,\\
    R^O_{\infty} & = & 0
  \end{array} \label{eq:bsde-f}
\end{equation}
A basic requirement for Eq.~\eqref{EqRemainderForceField} is that the source
term in~\eqref{EqRemainderForceField} and~\eqref{eq:bsde-f} is convergent in
the UV, i.e. \
\begin{equation}
  \label{EqIntegralConditionCoherence} \int_{a_0}^{\infty} \left\|| \bbE_{a_0}
  [\ring{\mathscr{L}}_c    \ring{f}_c (\varphi_c)] \right\|| \mathd c <
  \infty,
\end{equation}
for some scale $a_0$. If equation~\eqref{EqRemainderForceField} has indeed a
unique solution then $\ring{f}$ characterizes uniquely the force field $f$.

Assume now we are given an observable $O$ which has a germ $\ring{O}$. Under
proper assumptions on $\ring{O}$, and for a choice of $\ring{f}$ that ensure
an appropriate strong decay of $R^f_a$ as $a \rightarrow \infty$, the
coherence relation
\[ \int_b^{\infty} \left| \bbE_b [\mathscr{L}_c \ring{O}_c (\varphi_c)]
   \right| \hspace{0.17em} \mathd c < \infty \]
is a consequence of the relation
\[ \label{EqHatCoherence} \int_b^{\infty} \left| \bbE_b [\ring{\mathscr{L}_c} 
   \ring{O}_c (\varphi_c)] \right| \hspace{0.17em} \mathd c < \infty . \]
A family $\ring{O}$ that satisfies the estimate~\eqref{EqHatCoherence} is
called an \tmtextbf{\tmtextit{approximately coherent germ}}. In those terms,
Eq.~\eqref{EqIntegralConditionCoherence} states that $\ring{f}$ is an
approximately coherent germ for the force field $f$.

These considerations lead to the following strategy for
\tmtextit{constructing} the law of a random field $\varphi$. Associate to each
force germ $(\ring{f}_c)_c$ and each scale parameter $a_0 > 0$ the solution to
the FBSDEs
\begin{equation}
  \begin{array}{lll}
    \mathd \varphi^{a_0}_b & = & \dot{C}_b (\ring{f_b}  (\varphi^{a_0}_b) +
    R^{f, a_0}_b) \mathd b + \dot{C}_b^{1 / 2} \mathd W_b,\\
    R^{f, a_0}_b & = & \int_b^{a_0} \bbE_b  [R^{f, a_0}_c  \dot{C}_c \mathD
    \ring{f_c} (\varphi^{a_0}_c)] \mathd c\\
    &  & + \int_b^{a_0} \bbE_b [(\ring{\mathscr{L}}_c  \ring{f_c})
    (\varphi^{a_0}_c)] \mathd c
  \end{array} \label{EqCoupledSystemDynamicsRemainder}
\end{equation}
for $0 \leqslant b \leqslant a_0$ with mixed initial/final conditions
\begin{equation}
  \varphi_0^{a_0} = 0, \qquad R^{f, a_0}_{a_0} = 0.
\end{equation}
Now we need to find $(\ring{f_b} )_{b \geqslant 0}$ such that the coupled
system \eqref{EqCoupledSystemDynamicsRemainder} has a unique solution for all
$a_0$ and $a_0$-uniform estimates that entail sufficient compactness to pass
to the limit in \eqref{EqCoupledSystemDynamicsRemainder}.
In~{\cite{de_vecchi22}} this strategy has been used to construct some
subcritical Gibbs measures on Grassmann fields.

\subsection{Link with factorization algebras}

Definition~\ref{def:obs} gives rise only to a vector space structure on the
set of observables. Since the product of two martingales is, generally
speaking, not a martingale, we do not have a natural way of multiplying
observables. However for two observables $O^1, O^2$ that have some coherent
germs $\ring{O}^1, \ring{O}^2$ with disjoint supports, say they are at
distance strictly larger than $1 / a_0$, we can set for $a \geqslant a_0$
\[ \ring{O}^{(12)}_a \assign \ring{O}^{(1)}_a  \ring{O}^{(2)}_a . \]
Since $\tmop{Tr} \dot{C}_a^{1 / 2} \sigma_a^2 \dot{C}_a^{1 / 2} \mathD
\ring{O}^{(1)}_a \mathD \ring{O}^{(2)}_a = 0$, for $a \geqslant a_0$, because
of the support condition, one has
\[ \ring{\mathscr{L}_a}  \ring{O}^{(12)}_a = (\ring{\mathscr{L}_a} 
   \ring{O}^{(1)}_a) \ring{O}^{(2)}_a + \ring{O}^{(1)}_a (\ring{\mathscr{L}_a}
   \ring{O}^{(2)}_a) \]
therefore provided the germs are approximately coherent and we impose
sufficient decay on $\ring{\mathscr{L}_a}  \ring{O}^{(i)}_a$ and moderate
growth conditions on $\ring{O}^{(i)}_a$, it is possible to guarantee that
$\ring{\mathscr{L}_a}  \ring{O}^{(12)}_a$ is integrable and therefore
$\ring{O}^{(12)}_a$ is itself an approximately coherent germ which defines a
unique observable $O^{(12)} \backassign O^{(1)} \ast O^{(2)}$. This gives on a
subspace of observables a natural pre-factorization algebra structure, as
defined in
Costello--Gwilliam~{\cite{costelloFactorizationAlgebrasQuantum2017a,costelloFactorizationAlgebrasQuantum2021}}.

\section{Gauge theories}\label{SectionGauge}

In this section we assume that the field $\varphi$ is a connection on a
principal bundle over $\mathbb{R}^d$, with finite dimensional compact
structure group  G  and Lie algebra $\mathfrak{g}$. Recall that the space of
connections is affine with $\mathfrak{g}$-valued $1$-forms as underlying
vector space. The gauge group $\mathfrak{}$consists of  G -valued functions
and acts on connections by $g \cdummy \varphi = \tmop{Ad}_g \varphi - (\mathd
g) g^{- 1}$ and forms by $g \cdummy f = \tmop{Ad}_g f$. We assume moreover
that the force field is a gauge covariant function of the underlying field
$f_a = f_a (\varphi_a)$ (taking values in the space of $1$-forms) :
\[ g \cdot f_a (\psi) = f_a  (g \cdot \psi), \]
for every $g \in \mathfrak{G}$ and connection $\psi$. Given a connection
$\varphi$, denote by $h^{xy} (\varphi)$ the $\varphi$-holonomy along the
geodesic from $x$ to $y$. Recall that
\begin{equation}
  \label{EqActionOnHolonomy} h^{xy}  (g \cdot \varphi) = g (x) h^{xy}
  (\varphi) g (y)^{- 1} .
\end{equation}
Let $\chi_a (x, y)$ be a symmetric function of $(x, y)$. We define a map
$\dot{C}_a^{1 / 2} (\varphi)$ acting on $1$-forms by
\[ (\dot{C}_a^{1 / 2} (\varphi) \omega) (x) \assign \frac{1}{a^{1 / 2}}  \int
   \chi_a (x, y) \mathrm{Ad}_{h^{xy} (\varphi)} \omega (y) \mathd y. \]
The operator $\dot{C}_a^{1 / 2}$ is symmetric and $\dot{C}_a (\varphi) \assign
(\dot{C}_a^{1 / 2} (\varphi))^2$ is symmetric and non-negative. By
\eqref{EqActionOnHolonomy}, it is also gauge covariant:
\[ g \cdot (\dot{C}_a (\varphi) \omega) = \dot{C}_a  (g \cdot \varphi)  (g
   \cdot \omega) . \]
The tangent space of the gauge orbit at a given connection $\varphi$ is
spanned by the elements of the form $\mathd_{\varphi} h$, where
$\mathd_{\varphi}$ is the $\varphi$-covariant derivative and $h$ an arbitrary
(smooth enough) $\mathfrak{g}$-valued $0$-form,
\[ \mathd_{\varphi} h = \sum_i (\partial_i h + [\varphi_i, h]) \mathd x_i . \]
Two scale-dependent families of connections $(\varphi_a)_{a \geqslant 0}$ and
$(\varphi'_a)_{a \geqslant 0}$ are said to be \tmtextit{gauge equivalent} if
there exists a scale-dependent family of gauge transforms $(g_a)_{a \geqslant
0}$ such that $\varphi'_a = g_a \cdot \varphi_a$ for all $a \geqslant 0$. We
note that if $(\varphi_a)_{a \geqslant 0}$ is a solution of the
Wilson--It{\^o} equation \eqref{EqWIDE} with $\sigma_a \equiv 1$, and
$(g_a)_{a \geqslant 0}$ is an adapted process that is differentiable in $a$
and takes values in $C^1 (M, \mathfrak{G})$ then $\varphi^g_a \assign g_a
\cdot \varphi_a$ satisfies the equation
\begin{equation}
  \label{EqGaugedDynamics} \begin{array}{lll}
    \mathd \varphi^g_a & = & g_a \cdummy (\mathd \varphi_a) -
    \mathd_{\varphi_t^a} (\dot{g}_t g_t^{- 1} \mathd a)\\
    & = & (\dot{C}_a f_a (\varphi^g_a) - \mathd_{\varphi_a^g} (\dot{g}_a g^{-
    1}_a)) \mathd a\\
    &  & + \dot{C}_a^{1 / 2} (\varphi^g_a) \mathd W_a .
  \end{array}
\end{equation}
This is a particular case of a more general situation.

\begin{proposition}
  \label{LemGaugeAction}For any adapted process $(h_a)_{a \geqslant 0}$ with
  values in $C^1 (M, \mathfrak{g)}$ the solution $(\varphi^{(h)}_a)_{a
  \geqslant 0}$ to the equation
  \[ \mathd \varphi^{(h)}_a = ((\dot{C}_a f_a) (\varphi^{(h)}_a) +
     \mathd_{\varphi^{(h)}_a} h_a) \mathd a + \dot{C}_a^{1 / 2}
     (\varphi^{(h)}_a) \mathd W_a \]
  is gauge equivalent to the solution of a Wilson--It{\^o} equation
  \begin{equation}
    \mathd \varphi^{[h]}_a = (\dot{C}_a f_a) (\varphi^{[h]}_a) \mathd a +
    \dot{C}_a^{1 / 2} (\varphi^{[h]}_a) \mathd W^h_a \label{eq:g-transf}
  \end{equation}
  driven by another Brownian motion $W^h$. So the law of the gauge orbit of
  $(\varphi^{(0)}_a)_{a \geqslant 0}$ is well-defined.
\end{proposition}

The proof proceeds by solving the ordinary differential equation $\dot{g}_a
g_a^{- 1} = h_a$ and remarking that $\varphi^{[h]}_a \assign g_a \cdummy
\varphi^{(h) }_a$ solves~\eqref{eq:g-transf} with $\mathd W^h_a : = g_a \cdot
\mathd W_a$, which is a Brownian motion by It{\^o} isometry.

As Parisi--Wu stochastic quantization
scheme~{\cite{parisiPerturbationTheoryGauge1981}} this approach allows to
quantize gauge theories without path integrals and their associated ghosts or
BRS symmetries. Adding terms of the form $\mathd_{\varphi_a} h_a$ in the drift
allows us to perform gauge-fixing in analogy to the Zwanziger--DeTurck--Sadun
trick~{\cite{damgaardStochasticQuantization1988}}. Moreover it points to a
covariant formulation of the flow equation for the effective force $f_a$
\begin{equation}
  \partial_a f_a + f_a \dot{C}_a \mathD f_a + \frac{1}{2} \tmop{Tr} \dot{C}_a
  \mathD^2 f_a = 0 \label{eq:C-polchinski}
\end{equation}
with terminal condition $f_{\infty}$. This is similar to Polchinski
equation~\eqref{eq:polchinski} with the difference that the cutoff propagator
has been replaced by a covariant averaging operator.

\section{Conclusions}

We introduced Wilson--It{\^o} fields, a novel class of random fields described
by ``local interactions''. Our point of view is inscribed in the general idea
of \tmtextbf{\tmtextit{stochastic
quantization}}~{\cite{damgaardStochasticQuantization1988}} initiated by
Parisi--Wu~{\cite{parisiPerturbationTheoryGauge1981}}, which replaces the use
of functional integrals with stochastic partial differential equations for the
Euclidean fields. This topic has recently witnessed a renewed interest since
SPDEs can be used to study rigorously Euclidean fields: see e.g.~
{\cite{chandraLangevinDynamic2D2020,chandraStochasticQuantisationYangMillsHiggs2022,chevyrevInvariantMeasureUniversality2023}}
for non-Abelian gauge theories,
{\cite{bailleulPhiMeasuresCompact2023,bailleulUniquenessPhiMeasures2023}} for
$\Phi^4_3$ on Riemmanian manifolds and also~{\cite{GH21}} for a partial
verification of the Osterwalder--Schrader axioms in $\Phi^4_3$. However, at
variance with classical stochastic quantization methods, our approach does not
require to introduce fictious additional parameters and depends on a
physically relevant scale of observation, motivated by the Kadanoff--Wilson
picture of renormalization and by the mathematical framework for Euclidean QFT
via effective theories, see
e.g.~{\cite{costelloRenormalizationEffectiveField2011,salmhoferRenormalizationIntroduction2007,rychkovEPFLLecturesConformal2016}}.
Wilson--It{\^o} fields have a natural intrinsic definition, independent of
cutoff procedures and it seems interesting to pursue further their study, e.g.
investigate spatial Markov properties, reflection positivity,
uniqueness/non-uniqueness of solutions, numerical approximations, description
of scale-invariant and conformal invariant fields.

\end{document}